\documentclass[11pt,a4paper]{article}

\usepackage{inputenc}
\usepackage{amsmath}
\usepackage{bm}
\usepackage{bbold}
\usepackage{amsthm}
\usepackage{enumerate}

\usepackage{enumitem}

\usepackage{booktabs}
\usepackage{caption}

\usepackage{newfloat}
\DeclareFloatingEnvironment[name=Table]{specialtable}

\usepackage[top=1.0in, left=1.0in, right=1.0in, bottom=1.0in]{geometry}

\usepackage[hyphens]{url}

\usepackage{hyperref}
\usepackage{breakurl}

\title{Using pairwise comparisons to determine consumer preferences in hotel selection
\thanks{Mathematics, 2022, 10(5), 730, https://doi.org/10.3390/math10050730}
\thanks{This work was supported in part by the Russian Foundation for Basic Research (grant No. 20-010-00145).}}

\author{N. Krivulin\thanks{Faculty of Mathematics and Mechanics, St.~Petersburg State University, Universitetsky Ave.~28, St.~Petersburg, 198504, Russia; 
nkk@math.spbu.ru.}
\and
A. Prinkov\thanks{Faculty of Mathematics and Mechanics, St.~Petersburg State University, Universitetsky Ave.~28, St.~Petersburg, 198504, Russia; 
aprinkov@yahoo.com.}
\and 
I. Gladkikh\thanks{Graduate School of Management, St.~Petersburg State University, Universitetsky Ave.~28, St.~Petersburg, 198504, Russia; 
gladkikh@gsom.spbu.ru}
}

\date{}

\newtheorem{theorem}{Theorem}
\newtheorem{lemma}[theorem]{lemma}

\theoremstyle{definition}

\begin{document}

\maketitle

\begin{abstract}
We consider the problem of evaluating preferences for criteria used by university students when selecting a hotel for accommodation during a professional development program in a foreign country. Input data for analysis come from a survey of 202 respondents, who indicated their age, sex and whether they have previously visited the country. The criteria under evaluation are location, accommodation cost, typical guests, free breakfast, room amenities and courtesy of staff. The respondents assess the criteria both directly by providing estimates of absolute ratings and ranks, and indirectly by relative estimates using ratios of pairwise comparisons. To improve the accuracy of ratings derived from pairwise comparisons, we concurrently apply the principal eigenvector method, the geometric mean method and the method of log-Chebyshev approximation. Then, the results from the direct and indirect evaluation of ratings and ranks are examined together to analyze how the results from pairwise comparisons may differ from each other and from the results of direct assessment by respondents. We apply statistical techniques, such as estimation of means, standard deviations and correlations, to the vectors of ratings and ranks provided directly or indirectly by respondents, and then use the estimates to make accurate assessment of the criteria under study.
\\

\\

\end{abstract}

\section{Introduction}

Evaluation of preferences for alternatives based on their pairwise comparisons is a widely accepted approach in decision making, when direct assessment of the preferences is infeasible or impossible \cite{Thurstone1927Law,Saaty1990Analytic,Saaty2013Onthemeasurement,Gavalec2015Decision}. The approach uses the results of pairwise comparisons of alternatives on an appropriate scale, given in the form of a pairwise comparison matrix. Then, various computational methods can be applied to the matrix to make judgment on the preference of each alternatives by evaluating its individual rating (score, priority, weight) and ranking the alternatives according to the ratings.

The methods used to derive ratings from pairwise comparisons may exploit different computational techniques. These methods are mainly based on aggregation (summation) of columns in the pairwise comparison matrix to obtain a vector of ratings of alternatives, or on approximation of the pairwise comparison matrix by a symmetrically reciprocal (consistent) matrix that directly determines the vector of ratings (see, e.g., \cite{Saaty1984Comparison,Choo2004Common}). The key issue in deriving ratings from pairwise comparisons by different techniques is that these techniques may produce rather different or even contradictory results (see, e.g., \cite{Barzilai1987Consistent,Ishizaka2006How,Tran2013Pairwise,Mazurek2021Numerical}), which can significantly complicate or even make it impossible to unambiguously assess alternatives when making decisions in practice.

Among other problems arising from the issue indicated above, one has to recognize how much the results of different methods vary in practical problems, and how to make the right decision in case of possible inconsistency of the results. The comparison of results of assessment methods is discussed in many research studies \cite{Barzilai1997Deriving,Ishizaka2006How,Mazurek2021Numerical}, including statistical analysis of extensive data given by many pairwise comparison matrices. However, in most cases, the source data used in the analysis are obtained by simulation \cite{Ishizaka2006How,Mazurek2021Numerical} and thus are to be considered as an artificial input that may less adequately reflect the usual conditions of practice than authentic results of human evaluation. 

In practical situations where the methods used can give different results, a natural but reasonable way to evaluate alternatives more unambiguously is based on the simultaneous application of several methods to make a decision that concurrently considers all the results. If these results differ significantly, the choice of one of them as the basis for making a decision does not seem entirely justified. On the contrary, the closeness and stability of the results can serve as additional arguments in favor of choosing one of them as a solution that can be considered as close to optimal. 

The well-known solution approach is the principal eigenvector method \cite{Saaty1984Comparison,Saaty1990Analytic,Saaty2013Onthemeasurement}, which defines the vector of ratings as a weighted sum of the columns in the pairwise comparison matrix, where the weights are taken proportional to the entries of the vector of ratings itself. This approach leads to a solution calculated as the eigenvector of the pairwise comparison matrix, which corresponds to the maximal eigenvalue. The geometric mean method is another widely used approach, which solves the problem by approximating the pairwise comparison matrix by a consistent matrix in the Euclidean norm on the logarithmic scale~\mbox{\cite{Narasimhan1982Geometric,Crawford1985Note,Barzilai1987Consistent}.} It provides the solution vector by calculating the geometric mean of the entries in each row of the pairwise comparison matrix.

An approach to the derivation of ratings from pairwise comparison matrices, which applies the Chebyshev approximation on the logarithmic scale (the log-Chebyshev approximation), is developed in \cite{Krivulin2015Rating,Krivulin2016Using,Krivulin2018Methods,Krivulin2019Tropical}. The proposed solution technique is based on methods and results of tropical mathematics which is concerned with the study and application of algebraic systems with idempotent operations \cite{Baccelli1993Synchronization,Kolokoltsov1997Idempotent,Golan2003Semirings,Heidergott2006Maxplus}. The approximation problem is formulated and solved as an optimization problem in terms of tropical mathematics (a tropical optimization problem), which yields a complete analytical solution given in a parametric form ready for numerical computation and formal analysis.

The problem of evaluating alternatives from their pairwise comparisons often arises in studying consumer preferences in various markets including the market of hotel services. Specifically, for hotel suppliers, it is crucial to identify the most and least important criteria used by consumers in hotel selection. The problem of assessment criteria for hotel selection by various customer groups is studied in many researches (see, e.g., \cite{Tsaur1996Multiattribute,Dolnicar2003Which,Yavas2005Dimensions,Hsieh2008Service,Jones2011Factors,Chang2014Study,Yu2018Multicriteria,Tumer2019Assessment,Goral2020Prioritizing}). These studies mainly concentrate on investigation of preferences of business and (or) leisure travelers, who often constitute the largest segment of hotel guests and thus are of prime interest~\mbox{\cite{Dolnicar2003Which,Yavas2005Dimensions,Hsieh2008Service,Jones2011Factors,Goral2020Prioritizing},} but sometimes more specific customer groups such as university students are under consideration \cite{Tumer2019Assessment}.

As a source of information to understand how different factors (hotel attributes), such as location, accommodation cost, free breakfast or amenities, may affect the choice of a hotel, consumer surveys are widely used which can be administered online \cite{Jones2011Factors,Yu2018Multicriteria} or using a paper-and-pencil format \cite{Tsaur1996Multiattribute,Yavas2005Dimensions,Chang2014Study,Tumer2019Assessment,Goral2020Prioritizing}. A typical survey can ask respondents for absolute estimates that rate and rank the factors directly or for relative estimates in the form of pairwise comparison ratios for the factors.

The pairwise comparison data are then used to make a final assessment of factors by applying one of the methods of rating alternatives from pairwise comparisons. However, many studies rely on results obtained using only one method, which can lead to inaccurate or wrong conclusions because different methods may produce ambiguous results. To~improve the quality of the assessment, it seems necessary to verify the results of one method against the results of other available methods, which makes the development and experimental study of a combined approach to integrate several solution methods together with procedures of data analysis, a highly relevant issue. 

In this paper, we consider the problem of evaluating preferences for criteria used by university students when selecting a hotel for accommodation during a short-term professional development program in a foreign country. As input data for analysis, we are restricted by the results of a survey of 202 respondents who indicated their age, sex and whether they have previous experience of visiting the country. The criteria evaluated by respondents include location, accommodation cost, typical guests, free breakfast, room amenities and courtesy of staff. The respondents assess the criteria both directly by providing estimates of absolute ratings and ranks, and indirectly by relative estimates that are represented as ratios of pairwise comparisons.

The purpose of this applied study is twofold. First, it is to investigate the possibilities of combining different methods of rating and ranking criteria together with statistical procedures to provide more accurate assessment of the criteria, to describe an applied technique intended to improve the quality of assessment, and to demonstrate the implementation of this technique to the problem of evaluating criteria for hotel selection. The second objective is to explore by statistical procedures the results obtained by different methods from pairwise comparisons provided by respondents to gain additional insight and understanding about the similarity and difference between these results. 

As the key technique to improve the accuracy of ratings derived from the results of pairwise comparisons, we concurrently apply the methods of principal eigenvector, geometric mean and log-Chebyshev approximation. Then, the results obtained by the direct and indirect assessment of ratings and ranks are examined together to analyze how the results from pairwise comparisons may differ from each other and from the results of direct assessment. The analysis involves statistical techniques, such as estimation of means, standard deviations and correlations, applied to the vectors of ratings and ranks of criteria, which are provided directly or indirectly by respondents.

The paper is organized as follows. We start in Section~\ref{S-EAPC} with a brief overview of the methods used below to derive ratings from pairwise comparisons, including the method based on the log-Chebyshev approximation. Section~\ref{S-TMBS} presents the basic definitions, notation and results of tropical mathematics, which are then used to give direct formulas for calculating vectors of ratings using log-Chebyshev approximation. In Section~\ref{S-ECHS}, we formulate the problem of evaluating criteria for hotel selection, which motivates and illustrates the study, and describe the survey data used in the solution.

Section~\ref{S-SARR} offers results of statistical analysis of ratings and ranks derived from the survey data. First, we investigate whether characteristics of respondents (age, sex, etc.) may affect the degree of difference or similarity between vectors of ratings or ranks obtained for each respondent by different methods. Furthermore, statistical results of the analysis and comparison of the vectors of ratings for the criteria are given including estimates of means, standard deviations and correlations. We conclude the section with the results of comparison of rank vectors that order criteria according to their preferences. In Section~\ref{S-D}, we present some discussions, and in Section~\ref{S-C} offer concluding remarks.

\section{Evaluation of Alternatives from Pairwise Comparisons}
\label{S-EAPC}

Suppose there are $n$ alternatives $\mathcal{A}_{1},\ldots,\mathcal{A}_{n}$ for making decisions, which are compared in pairs. The results of the comparisons are represented in the form of a pairwise comparison matrix $\bm{A}=(a_{ij})$ of dimension $n\times n$, where the element $a_{ij}>0$ indicates by how many times the alternative $\mathcal{A}_{i}$ is more preferable than $\mathcal{A}_{j}$. The problem of pairwise comparisons is to find a vector $\bm{x}=(x_{i})$ of absolute ratings (scores, priorities, weights) of alternatives on the basis of the relative estimates given by pairwise comparisons.

Consider a pairwise comparison matrix $\bm{A}$ and note that its elements satisfy the condition $a_{ji}=1/a_{ij}$ (and hence $a_{ii}=1$) for all $i,j=1,\ldots,n$, and therefore the matrix $\bm{A}$ is symmetrically reciprocal. Furthermore, if the transitive property in the form of equality $a_{ij}=a_{ik}a_{kj}$ is fulfilled for all $i,j,k=1,\ldots,n$, then the matrix $\bm{A}$ is called consistent.

The elements of a consistent pairwise comparison matrix $\bm{A}$ are represented as $a_{ij}=x_{i}/x_{j}$ for all $i,j=1,\ldots,n$, where $x_{i}$ are the components of some positive vector $\bm{x}$ determined up to a positive factor. Moreover, it directly follows from the equality $a_{ij}=x_{i}/x_{j}$ that $\bm{x}=(x_{i})$ is a solution of the pairwise comparison problem of interest.

In practical applications, the matrices obtained by pairwise comparisons of alternatives are usually not consistent. In this case, a problem arises of finding a consistent matrix that is close to (approximates) the original pairwise comparison matrix. Any vector that determines this consistent matrix is then taken as the vector of absolute ratings. To find an appropriate consistent matrix of pairwise comparisons various heuristic procedures and approximation methods are used (see, e.g., \cite{Saaty1984Comparison,Saaty1990Analytic,Choo2004Common}).

Heuristic methods for solving the problem usually apply techniques of aggregation (summation) of the columns in the pairwise comparison matrix. The columns are added up with coefficients (weights) that are dependent on the technique. The resulting vector, which generates some consistent matrix, serves as a solution to the problem. As examples of the heuristic approach, one can consider the weighted column sum method and the principal eigenvector method.

The most common in practice is the method of principal eigenvector proposed by T.~Saaty in the 1970s \cite{Saaty1977Scaling,Saaty1984Comparison,Saaty1990Analytic,Saaty2013Onthemeasurement}. The method assumes that the weights of columns for aggregation are set proportional to the components of the solution vector of absolute ratings. This leads to a solution of the problem that takes the form of the principal eigenvector $\bm{x}$ of the pairwise comparison matrix $\bm{A}$, which corresponds to the maximum eigenvalue $\lambda$ of the matrix and satisfies the equality
\begin{equation*}
\bm{A}\bm{x}
=
\lambda\bm{x}.
\end{equation*}

The approximation methods solve a problem of minimizing an error in the approximation of the pairwise comparison matrix $\bm{A}=(a_{ij})$ by a consistent matrix $\bm{X}=(x_{i}/x_{j})$. The application of these methods offers a mathematically justified approach, which however can result in very hard optimization problems. The complexity of the solution essentially depends on the metric and scale used to estimate the approximation error. 

To measure the error on the standard linear scale, different metrics may be used including the Euclidean, Manhattan and Chebyshev distance functions. However, minimizing the errors in linear scale usually leads to complex nonlinear multiextremal optimization problems that are difficult to solve \cite{Saaty1984Comparison,Chu1998Ontheoptimal}, and hence is not common in practice.

If the approximation error is evaluated in logarithmic scale, then the solution may become less complicated and can sometimes be obtained in analytical form \cite{Narasimhan1982Geometric,Crawford1985Note,Barzilai1987Consistent}. Note that when transforming from linear scale to logarithmic scale, the variation between entries in a pairwise comparison matrix $\bm{A}$ decreases, which reduces the impact on the overall error of large values (for example, $a_{ij}=9$) to the detriment of small ones ($a_{ji}=1/9$). Therefore, the use of the logarithmic scale in the approximation of the pairwise comparison matrices seems to be quite reasonable. 

Consider the application of the $\log$-Euclidean metric (the Euclidean metric on the logarithmic scale), in which the distance between the matrices $\bm{A}=(a_{ij})$ and $\bm{X}=(x_{i}/x_{j})$ is determined using logarithm in a base greater than one by the formula
\begin{equation*}
l_{2}(\bm{A},\bm{X})
=
\Bigg(
\sum_{1\leq i,j\leq n}
\left(
\log a_{ij}-\log\frac{x_{i}}{x_{j}}
\right)^{2}
\Bigg)^{1/2}.
\end{equation*}

The minimum distance is immediately found by calculating the derivatives of the squared distance $l_{2}^{2}(\bm{A},\bm{X})$ with respect to all $x_{i}$, equating the derivatives to zero, and solving the obtained equations. For a symmetrically reciprocal matrix $\bm{A}$, the result is a direct solution of the approximation problem, in which the components of the vector $\bm{x}=(x_{1},\ldots,x_{n})^{T}$ that determines the matrix $\bm{X}$ are given in the parametric form
\begin{equation*}
x_{i}
=
\Bigg(
\prod_{j=1}^{n}a_{ij}
\Bigg)^{1/n}
u,
\qquad
u>0,
\qquad
i=1,\ldots,n.
\end{equation*}

The obtained vector $\bm{x}$ (usually normalized in rectangular metric, i.e., with respect to the sum of components) is called the solution of the problem of pairwise comparisons by the method of geometric mean \cite{Narasimhan1982Geometric,Crawford1985Note,Barzilai1987Consistent}.

The distance between the matrices $\bm{A}$ and $\bm{X}$ in the Chebyshev metric in logarithmic scale (the log-Chebyshev distance) is defined as
\begin{equation*}
l_{\infty}(\bm{A},\bm{X})
=
\max_{1\leq i,j\leq n}
\left|\log a_{ij}-\log\frac{x_{i}}{x_{j}}\right|.
\end{equation*}

The problem of approximating a pairwise comparison matrix in the log-Chebyshev metric can be reduced to the following problem of minimizing a function without logarithm (see, e.g., \cite{Krivulin2019Tropical,Krivulin2020Using}):
\begin{equation}
\begin{aligned}
&
\min_{\bm{x}>\bm{0}}
&&
\max_{1\leq i,j\leq n}
\frac{a_{ij}x_{j}}{x_{i}},
\end{aligned}
\label{P-minxmaxaijxjxi}
\end{equation}
which is equivalent to minimizing the maximum relative error \cite{Elsner2004Maxalgebra,Krivulin2020Using}, given by
\begin{equation*}
\max_{1\leq i,j\leq n}
\left|\frac{a_{ij}-x_{i}/x_{j}}{a_{ij}}\right|.
\end{equation*}

In contrast to both methods of the principal eigenvalue and geometric mean, which always lead to a single solution vector (up to a positive factor), the solution based on the log-Chebyshev approximation may be nonunique.

The existence of multiple solutions instead of a single one, on the one hand, can make it difficult to choose the most preferable alternative in practice. On the other hand, the availability of a set of different solutions to the problem extends the possibility of making optimal decisions, for example by taking into account additional constraints to narrow the solutions. In addition, due to the rather approximate nature of the model of pairwise comparisons for which the inconsistency of pairwise judgments is typical, the existence of several solutions to the problem seems to be quite natural.

Suppose that the approximation procedure results in a set $\mathcal{S}$ of optimal vectors, rather than in a single optimal vector of ratings. As the ``best'' and ``worst'' solutions to the problem, we can take those two vectors from $\mathcal{S}$ that best and worst differentiate the alternatives with highest and lowest ratings \cite{Krivulin2018Methods,Krivulin2019Tropical}. 

Consider the ratio between the maximum and minimum elements of the vector $\bm{x}=(x_{i})$, which is called the Hilbert (interval, range) seminorm, and given by
\begin{equation*}
\max_{1\leq i\leq n}x_{i}\Big/\min_{1\leq j\leq n}x_{j}
=
\max_{1\leq i\leq n}x_{i}\times\max_{1\leq j\leq n}x_{j}^{-1}.
\end{equation*}

The solution vector with the maximum Hilbert seminorm is taken as the best differentiating solution, whereas the vector with the minimum seminorm is as the worst differentiating solution. These vectors are obtained by solving the following problems:
\begin{gather}
\begin{aligned}
&
\max_{\bm{x}\in\mathcal{S}}
&&
\max_{1\leq i\leq n}x_{i}\times\max_{1\leq j\leq n}x_{j}^{-1},
\end{aligned}
\label{P-maxxSmaxximaxxj}
\\
\begin{aligned}
&
\min_{\bm{x}\in\mathcal{S}}
&&
\max_{1\leq i\leq n}x_{i}\times\max_{1\leq j\leq n}x_{j}^{-1}.
\end{aligned}
\label{P-minxSmaxximaxxj}
\end{gather}

Problems \eqref{P-minxmaxaijxjxi}--\eqref{P-minxSmaxximaxxj}, which arise from the log-Chebyshev approximation in the pairwise comparison problem, can be directly solved in analytical form, using methods and results of tropical mathematics (see, e.g., \cite{Krivulin2015Rating,Krivulin2016Using,Krivulin2018Methods,Krivulin2019Tropical} for further details).

\section{Tropical Mathematics Based Solutions}
\label{S-TMBS}

The purpose of this section is to provide a brief overview of the basic concepts, definitions and notation of tropical mathematics, which are used to describe the solution of the pairwise comparison problem. Further information on the theory, methods and applications of tropical mathematics can be found, for example, in \cite{Baccelli1993Synchronization,Cuninghamegreen1994Minimax,Kolokoltsov1997Idempotent,Golan2003Semirings,Heidergott2006Maxplus}.

\subsection{Elements of Tropical Mathematics}

Tropical (idempotent) mathematics deals with the theory and applications of algebraic systems with idempotent operations. An operation is idempotent if its application to arguments of the same value yields this value as the result. For example, taking the maximum is an idempotent operation since $\max\{x,x\}=x$, whereas the arithmetic addition is not: $x+x=2x$.

Tropical optimization focuses on optimization problems that are formulated and solved in terms of tropical mathematics. Models and methods of tropical optimization make it possible to find new solutions for classical and newly posed problems. Many problems can be solved directly in explicit analytic form; for other problems, only algorithmic techniques that offer solutions in numerical form are known. Applications of tropical optimization include various problems in project scheduling, location analysis, decision making and other areas.

An example of an algebraic system with an idempotent operation is the max-algebra defined on the set of non-negative real numbers $\mathbb{R}_{+}=\{x\in\mathbb{R}|x\geq0\}$. It is closed under addition denoted by $\oplus$ and defined for all $x,y\in\mathbb{R}_{+}$ as maximum: $x\oplus y=\max\{x,y\}$, and multiplication defined as usual. The neutral elements with respect to addition and multiplication coincide with the arithmetic zero $0$ and one $1$.

The tropical addition $\oplus$ is not invertible (opposite numbers do not exist), and hence a subtraction operation is undefined in max-algebra. The notion and notation of inverse elements with respect to multiplication, and exponents have the usual sense.

Vector and matrix operations are performed according to standard rules, where the arithmetic addition $+$ is replaced by $\oplus$. Specifically, multiplication of a vector or matrix by a scalar has the same result as in the standard arithmetic. In what follows, all vectors are considered column vectors until otherwise indicated. The zero vector denoted by $\bm{0}$, positive vector, zero matrix and identity matrix denoted by $\bm{I}$ are defined as usual.

For any nonzero column vector $\bm{x}=(x_{j})$, its multiplicative conjugate transpose is the row vector $\bm{x}^{-}=(x_{j}^{-})$, where $x_{j}^{-}=x_{j}^{-1}$ if $x_{j}\ne0$, and $x_{j}^{-}=0$ otherwise. For the vector of all ones, which is denoted as $\bm{1}$, the conjugate transpose is given by $\bm{1}^{-}=\bm{1}^{T}$. Multiplicative conjugate transposition of a nonzero matrix $\bm{A}=(a_{ij})$ results in the transposed matrix $\bm{A}^{-}=(a_{ij}^{-})$, where $a_{ij}^{-}=a_{ji}^{-1}$ if $a_{ji}\ne0$, and $a_{ij}^{-}=0$ otherwise.

For vectors $\bm{a}_{1},\ldots,\bm{a}_{n}$, its linear combination with nonnegative coefficients $x_{1},\ldots,x_{n}$ is given by $x_{1}\bm{a}_{1}\oplus\cdots\oplus x_{n}\bm{a}_{n}$. A vector $\bm{b}$ linearly depends on vectors $\bm{a}_{1},\ldots,\bm{a}_{n}$ if there exist numbers $x_{1},\ldots,x_{n}$ such that $\bm{b}=x_{1}\bm{a}_{1}\oplus\cdots\oplus x_{n}\bm{a}_{n}$. The collinearity of two vectors has the usual sense: vectors $\bm{b}$ is collinear with $\bm{a}$ if $\bm{b}=x\bm{a}$ for some $x$.

The set of all linear combinations $x_{1}\bm{a}_{1}\oplus\cdots\oplus x_{n}\bm{a}_{n}$ of a system of vectors $\bm{a}_{1},\ldots,\bm{a}_{n}$ forms a tropical linear space. Any vector $\bm{y}$ of this space is expressed as the (tropical) product of the matrix $\bm{A}=(\bm{a}_{1},\ldots,\bm{a}_{n})$, which consists of the vectors in the system taken as columns, and some vector $\bm{x}=(x_{1},\ldots,x_{n})^{T}$, in the form $\bm{y}=\bm{A}\bm{x}$.

For a square matrix $\bm{A}=(a_{ij})$ of order $n$, the nonnegative integer powers indicate repeated tropical multiplication of the matrix with itself, and are defined as $\bm{A}^{0}=\bm{I}$, $\bm{A}^{p}=\bm{A}^{p-1}\bm{A}=\bm{A}\bm{A}^{p-1}$ for all integer $p>0$. The trace of the matrix is given by
\begin{equation*}
\mathop\mathrm{tr}\bm{A}=a_{11}\oplus\cdots\oplus a_{nn}
=
\bigoplus_{i=1}^{n}a_{ii}.
\end{equation*}

The spectral radius of the matrix $\bm{A}$ is calculated by the formula
\begin{equation*}
\lambda
=
\mathop\mathrm{tr}\bm{A}\oplus\cdots\oplus\mathop\mathrm{tr}\nolimits^{1/n}(\bm{A}^{n})
=
\bigoplus_{k=1}^{n}{\mathop\mathrm{tr}}^{1/k}(\bm{A}^{k}).
\end{equation*}

On condition that $\lambda\leq1$, the Kleene operator (the Kleene star) is defined to map the matrix $\bm{A}$ onto the matrix
\begin{equation*}
\bm{A}^{\ast}
=
\bm{I}\oplus\bm{A}\oplus\cdots\oplus\bm{A}^{n-1}
=
\bigoplus_{k=0}^{n-1}\bm{A}^{k}.
\end{equation*}

\subsection{Preliminary Results}

We start with a formal criterion of linear dependence of vectors. In order to check whether a vector $\bm{b}$ is dependent on the system of vectors $\bm{a}_{1},\ldots,\bm{a}_{n}$, the following result can be used \cite{Krivulin2006Solution} (see, also \cite{Cuninghamegreen1994Minimax}).
\begin{lemma}
\label{L-AbAbeq1}
Denote by $\bm{A}=(\bm{a}_{1},\ldots,\bm{a}_{n})$ a matrix that consists of vectors $\bm{a}_{1},\ldots,\bm{a}_{n}$ as columns. Then, a vector $\bm{b}$ is linearly dependent on $\bm{a}_{1},\ldots,\bm{a}_{n}$ if and only if the following equality holds:
\begin{equation*}
(\bm{A}(\bm{b}^{-}\bm{A})^{-})^{-}\bm{b}
=
1.
\end{equation*}
\end{lemma}

We now consider some results in tropical optimization that provide the basis to solve pairwise comparison problems using log-Chebyshev approximation. Suppose that, given an $(n\times n)$-matrix $\bm{A}$, we need to find positive $n$-vectors $\bm{x}$ to solve the problem
\begin{equation}
\begin{aligned}
&
\min_{\bm{x}>\bm{0}}
&&
\bm{x}^{-}\bm{A}\bm{x}.
\end{aligned}
\label{P-minxxAx}
\end{equation}

A complete solution of the problem is obtained as follows (see, e.g., \cite{Krivulin2015Extremal}).
\begin{lemma}
\label{L-minxxAx}
Let $\bm{A}$ be a matrix with spectral radius $\lambda>0$. Then, the minimum in problem \eqref{P-minxxAx} is equal to $\lambda$, and all positive solutions are given in the parametric form
\begin{equation*}
\bm{x}
=
(\lambda^{-1}\bm{A})^{\ast}
\bm{u},
\qquad
\bm{u}>\bm{0}.
\end{equation*}
\end{lemma}

Let $\bm{B}=(\bm{b}_{j})$ be a positive $(n\times m)$-matrix with columns $\bm{b}_{j}=(b_{ij})$, and the problem is to obtain positive $n$-vectors $\bm{x}$ that provide the maximum
\begin{equation}
\begin{aligned}
&
\max_{\bm{x}>\bm{0}}
&&
\bm{1}^{T}\bm{x}\bm{x}^{-}\bm{1},
\\
&&&
\bm{x}
=
\bm{B}
\bm{u},
\quad
\bm{u}>\bm{0}.
\end{aligned}
\label{P-maxx1xx1-xeqBu-uge0}
\end{equation}

The next result offers a direct solution to the problem \cite{Krivulin2016Maximization,Krivulin2018Methods,Krivulin2019Tropical}.
\begin{lemma}
\label{L-maxx1xx1-xeqBu-uge0}
Denote by $\bm{B}_{lk}$ the matrix obtained from $\bm{B}$ by setting to zero all elements except $b_{lk}$ for some indices $k$ and $l$. Then, the maximum in problem \eqref{P-maxx1xx1-xeqBu-uge0} is equal to $\Delta=\bm{1}^{T}\bm{B}\bm{B}^{-}\bm{1}$, and all positive solutions are given in the parametric form
\begin{equation*}
\bm{x}
=
\bm{B}(\bm{I}\oplus\bm{B}_{lk}^{-}\bm{B})
\bm{u},
\qquad
\bm{u}>\bm{0},
\end{equation*}
where the indices $k$ and $l$ are selected by the conditions
\begin{equation*}
k
=
\arg\max_{j}\bm{1}^{T}\bm{b}_{j}\bm{b}_{j}^{-}\bm{1},
\qquad
l
=
\arg\max_{i}b_{ik}^{-1}.
\end{equation*}
\end{lemma}

Finally, suppose that, given an $(n\times n)$-matrix $\bm{A}$ with spectral radius $\lambda>0$, we need to find positive $n$-vectors $\bm{x}$ that attain the minimum
\begin{equation}
\begin{aligned}
&
\min_{\bm{x}>\bm{0}}
&&
\bm{1}^{T}\bm{x}\bm{x}^{-}\bm{1},
\\
&&&
\bm{x}
=
(\lambda^{-1}\bm{A})^{\ast}
\bm{u},
\quad
\bm{u}>\bm{0}.
\end{aligned}
\label{P-minx1xx1-xeqlambda1Aast-uge0}
\end{equation}

A solution to this problem is found as follows \cite{Krivulin2018Methods,Krivulin2019Tropical}.
\begin{lemma}
\label{L-minx1xx1-xeqlambda1Aast-uge0}
Let $\bm{A}$ be a matrix with spectral radius $\lambda>0$. Then, the minimum in problem \eqref{P-minx1xx1-xeqlambda1Aast-uge0} is equal to $\delta=\bm{1}^{T}(\lambda^{-1}\bm{A})^{\ast}\bm{1}$, and all positive solutions are given in the parametric form
\begin{equation*}
\bm{x}
=
(\delta^{-1}\bm{1}\bm{1}^{T}\oplus\lambda^{-1}\bm{A})^{\ast}
\bm{u},
\qquad
\bm{u}>\bm{0}.
\end{equation*}
\end{lemma}

Finally, note that the parametric form of solutions offered by Lemmas~\ref{L-minxxAx}--\ref{L-minx1xx1-xeqlambda1Aast-uge0} defines the set of solutions to problems \eqref{P-minxxAx}--\eqref{P-minx1xx1-xeqlambda1Aast-uge0} as a linear span of columns of corresponding matrices that generate the solutions.

\subsection{Application to Pairwise Comparison Problem}

Consider the problem of evaluating ratings of alternatives from pairwise comparison, which is solved using the log-Chebyshev approximation of a pairwise comparison matrix. In this case, we need to obtain a solution of problem \eqref{P-minxmaxaijxjxi}, and then, if the result is not unique, solutions of problems \eqref{P-maxxSmaxximaxxj} and \eqref{P-minxSmaxximaxxj}. Below, we use results in \cite{Krivulin2015Rating,Krivulin2016Using,Krivulin2018Methods,Krivulin2019Tropical} to demonstrate how these problems can be analytically solved in explicit form. 

We combine the solutions of the above problems into a procedure to handle the pairwise comparison problem under consideration. The procedure includes the following main steps: (i) finding the set of all solution vectors that represent absolute ratings of alternatives; (ii) checking whether the vectors determine a unique solution of the problem; and, in the case of a nonunique solution, (iii) selecting vectors which best and worst differentiate alternatives with the highest and lowest ratings (the best and worth differentiating~solutions).

\subsubsection{Finding All Solution Vectors}

First, we observe that the calculation of all solutions to the pairwise comparison problem with a matrix $\bm{A}$ by using log-Chebyshev approximation requires solving problem~\eqref{P-minxmaxaijxjxi}. In the tropical algebra setting, the objective function in this problem takes the vector form $\bm{x}^{-}\bm{A}\bm{x}$, whereas the problem itself is written as \eqref{P-minxxAx} and can be solved by Lemma~\ref{L-minxxAx}. The application of the lemma starts with the evaluation of the spectral radius for the matrix $\bm{A}$, given by
\begin{equation}
\lambda
=
\bigoplus_{k=1}^{n}{\mathop\mathrm{tr}}^{1/k}(\bm{A}^{k}).
\label{E-lambdaeqtr1kAk}
\end{equation}

Furthermore, we construct the matrix $\lambda^{-1}\bm{A}$ and calculate its Kleene star matrix
\begin{equation}
\bm{B}
=
(\lambda^{-1}\bm{A})^{\ast}
=
\bigoplus_{k=0}^{n-1}
(\lambda^{-1}\bm{A})^{k}.
\label{E-Beqlambda1Aasteqlambda1Ak}
\end{equation}

The set of all solutions of the problem is described using a vector of parameters $\bm{u}$ as
\begin{equation}
\bm{x}
=
\bm{B}\bm{u}.
\qquad
\bm{u}>\bm{0}
\label{E-xeqlambdaAastu-ugt0}
\end{equation}

This result shows that these solutions form a set of linear combinations of columns (are generated by columns) in the Kleene matrix $\bm{B}=(\lambda^{-1}\bm{A})^{\ast}$.

\subsubsection{Checking Uniqueness of Solution}

In the framework of the pairwise comparison problem, the solution obtained is considered unique if each column in the matrix $\bm{B}$ generates the same space of vectors, which happens when all columns are collinear. To check the uniqueness of the solution, we apply a refinement procedure based on Lemma~\ref{L-AbAbeq1}. The procedure examines the columns of the Kleene matrix $\bm{B}=(\bm{b}_{j})$ sequentially one by one and deletes a column if it is linearly dependent on others. For each column $\bm{b}_{j}$, we form a matrix $\bm{B}_{(j)}$ by taking all other columns of $\bm{B}$, which have not yet been deleted, and then verify the condition
\begin{equation*}
(\bm{B}_{(j)}(\bm{b}_{j}^{-}\bm{B}_{(j)})^{-})^{-}\bm{b}_{j}
=
1.
\end{equation*}

If this condition holds, the column $\bm{b}_{j}$ is deleted from the matrix $\bm{B}$, and otherwise retained. The solution is unique if the procedure results in a matrix with a single column that is taken as the solution. In the case when a matrix with several linearly independent columns is obtained to generate a set of different solutions, we turn to finding the best and worst differentiating vectors of ratings among them.

\subsubsection{Best Differentiating Solution}

Suppose that a solution to the log-Chebyshev approximation problem is obtained in the form $\bm{x}=\bm{B}\bm{u}$, where the matrix $\bm{B}$ has linearly independent columns. To select the best differentiating solution, we need to solve problem \eqref{P-maxxSmaxximaxxj}. In the framework of tropical algebra, this problem has the objective function written in the form $\bm{1}^{T}\bm{x}\bm{x}^{-}\bm{1}$, and can be represented as problem \eqref{P-maxx1xx1-xeqBu-uge0}, which is solved using Lemma~\ref{L-maxx1xx1-xeqBu-uge0}.

According to Lemma~\ref{L-maxx1xx1-xeqBu-uge0}, we need to find indices $k$ and $l$ that satisfy the conditions
\begin{equation*}
k
=
\arg\max_{j}\bm{1}^{T}\bm{b}_{j}\bm{b}_{j}^{-}\bm{1},
\qquad
l
=
\arg\max_{i}b_{ik}^{-1}.
\end{equation*}

Then, we construct the matrix $\bm{B}_{lk}$ by setting to zero all elements in the matrix $\bm{B}$ except for the element $b_{lk}$. The best differentiating solution is then given using a vector of parameters $\bm{u}_{1}$ by
\begin{equation*}
\bm{x}_{1}
=
\bm{B}(\bm{I}\oplus\bm{B}_{lk}^{-}\bm{B})
\bm{u}_{1},
\qquad
\bm{u}_{1}>\bm{0}.
\end{equation*}

\subsubsection{Worst Differentiating Solution}

We obtain the worst differentiating solution by solving problem \eqref{P-minxSmaxximaxxj}, which in terms of tropical algebra, takes the form of \eqref{P-minx1xx1-xeqlambda1Aast-uge0}. The solution of the latter problem by applying Lemma~\ref{L-minx1xx1-xeqlambda1Aast-uge0} involves evaluating
\begin{equation*}
\delta
=
\bm{1}^{T}(\lambda^{-1}\bm{A})^{\ast}\bm{1}.
\end{equation*}

The worst differentiating solution is given in parametric form as
\begin{equation*}
\bm{x}_{2}
=
(\delta^{-1}\bm{1}\bm{1}^{T}\oplus\lambda^{-1}\bm{A})^{\ast}\bm{u}_{2},
\qquad
\bm{u}_{2}
>
\bm{0}.
\end{equation*}

\section{Evaluation of Criteria for Hotel Selection}
\label{S-ECHS}

In this section, we consider the application of the above-discussed techniques of ratings alternatives from pairwise comparisons to the problem of evaluating consumer preferences in hotel selection. The problem is to identify priorities among criteria used by students in choosing hotels in a foreign country. The research is based on data obtained from 202 students of both sexes aged 17 to 26 years, who are to take a short-term professional development program at an educational institution in the country. Some students have visited this country and may have some experience of living there.

The students are provided with a grant that covers the travel expenses, tuition payments and food expenses during the working day. At the same time, the students must pay for accommodation on their own and choose a hotel by applying information available in the Internet on booking portals and hotel websites as well as their previous experience on visiting the country.

There are six criteria that are used when choosing a hotel: $\mathcal{C}_{1}$---location, $\mathcal{C}_{2}$---accom-modation cost, $\mathcal{C}_{3}$---social environment (typical guests who usually stay at the hotel), $\mathcal{C}_{4}$---free breakfast, $\mathcal{C}_{5}$---amenities (room equipment), $\mathcal{C}_{6}$---courtesy of hotel staff. The aim of the research is to evaluate ratings (scores, priorities, weights) and then obtain ranks of the criteria, based on the results of the survey of the students as respondents.

During the survey process, the respondent directly assesses the degree of importance (the rating) of each criterion, determines the ranks of the criteria, and compares the criteria in pairs. For the direct rating of criteria, the following absolute scale is used: $1/5$---not important, $2/5$---important, $3/5$---fairly important, $4/5$---important, $1$---very important. For the pairwise comparison of criteria, the following comparative scale is used: $1$---equal importance, $2$---moderate importance, $3$---strong importance, $4$---very strong importance; $5$---extreme importance of one criterion over another.

The survey results together with the parameters of respondents (sex, age, etc.) are summarized in a table. For each respondent $i$, the survey results include
\begin{description}[noitemsep, left=\parindent]
\item
[$\bm{s}_{\textbf{R}}(i)$,]
the vector of direct ratings (scores) of criteria by the respondent, 
\item
[$\bm{r}_{\textbf{R}}(i)$,]
the vector of direct ranks of criteria by the respondent, 
\item
[$\bm{C}_{\textbf{R}}(i)$,]
the matrix of pairwise comparisons of criteria, provided by the respondent. 
\end{description}

Below, the survey results for three respondents with $i=1,2,3$ are given as examples.

Respondent 1, male, age 24, visited the country before, provides the following estimates:
\begin{gather*}
\bm{s}_{\textbf{R}}(1)
=
\left(
\begin{array}{c}
3/5
\\
1
\\
1/5
\\
2/5
\\
1
\\
1/5
\end{array}
\right)
=
\left(
\begin{array}{c}
0.60
\\
1.00
\\
0.20
\\
0.40
\\
1.00
\\
0.20
\end{array}
\right),
\qquad
\bm{r}_{\textbf{R}}(1)
=
\left(
\begin{array}{c}
3
\\
2
\\
5
\\
4
\\
1
\\
6
\end{array}
\right),
\\
\bm{C}_{\textbf{R}}(1)
=
\left(
\begin{array}{cccccc}
 1 & 1/4 & 5 & 4 & 1/5 & 4
\\	
 4 & 1 & 5 & 5 & 1/3 & 5
\\	
1/5 & 1/5 & 1 & 1/3 & 1/5 & 2
\\
1/4 & 1/5 & 3 & 1 & 1/5 & 3
\\
 5 & 3 & 5 & 5 & 1 & 5
\\	
1/4 & 1/5 & 1/2 & 1/3 & 1/5 & 1
\end{array}
\right);
\end{gather*}

Respondent 2, male, age 19, never visited the country before,
\begin{gather*}
\bm{s}_{\textbf{R}}(2)
=
\left(
\begin{array}{c}
1
\\
1
\\
3/4
\\
1/2
\\
3/4
\\
1/2
\end{array}
\right)
=
\left(
\begin{array}{c}
1.00
\\
1.00
\\
0.75
\\
0.50
\\
0.75
\\
0.50
\end{array}
\right),
\qquad
\bm{r}_{\textbf{R}}(2)
=
\left(
\begin{array}{c}
2
\\
1
\\
3
\\
4
\\
5
\\
6
\end{array}
\right),
\\
\bm{C}_{\textbf{R}}(2)
=
\left(
\begin{array}{cccccc}
 1 & 1 & 3 & 3 & 2 & 4
\\
 1 & 1 & 4 & 2 & 4 & 5
\\
1/3 & 1/4 & 1 & 3 & 1 & 4
\\
1/3 & 1/2 & 1/3 & 1 & 1 & 3
\\
1/2 & 1/4 & 1 & 1 & 1 & 3
\\
1/4 & 1/5 & 1/4 & 1/3 & 1/3 & 1
\end{array}
\right);
\end{gather*}

Respondent 3, female, age 21, never visited the country before,
\begin{gather*}
\bm{s}_{\textbf{R}}(3)
=
\left(
\begin{array}{c}
4/5
\\
1
\\
1/5
\\
3/5
\\
4/5
\\
3/5
\end{array}
\right)
=
\left(
\begin{array}{c}
0.80
\\
1.00
\\
0.20
\\
0.60
\\
0.80
\\
0.60
\end{array}
\right),
\qquad
\bm{r}_{\textbf{R}}(3)
=
\left(
\begin{array}{c}
2
\\
1
\\
6
\\
4
\\
3
\\
5
\end{array}
\right),
\\
\bm{C}_{\textbf{R}}(3)
=
\left(
\begin{array}{cccccc}
 1 & 1/4 & 5 & 4 & 1/3 & 3
\\
 4 & 1 & 5 & 5 & 3 & 5
\\
1/5 & 1/5 & 1 & 1/3 & 1/5 & 1/3
\\
1/4 & 1/5 & 3 & 1 & 1/4 & 1
\\
 3 & 1/3 & 5 & 4 & 1 & 5
\\
1/3 & 1/5 & 3 & 1 & 1/5 & 1
\end{array}
\right).
\end{gather*}

Given the pairwise comparison matrix obtained from each respondent, the absolute ratings of criteria are found based on three computational methods. We obtain vectors of ratings by application of the principal eigenvector method and the method of geometric mean. Furthermore, we solve the problem by using the log-Chebyshev approximation in the framework of tropical algebra, which yields two vectors that are best and worst differentiate the criteria with the highest and lowest ratings (the best and worst log-Chebyshev approximation vectors). The vectors of ratings obtained are then used to rank criteria according to the values of ratings.

As a result, for each respondent $i$, the following vectors of ratings normalized by dividing by the maximum element, and rank vectors are calculated:
\begin{description}[noitemsep, left=\parindent]
\item
[$\bm{s}_{\textbf{R}}(i)$,]
the direct ratings (scores) by the respondent (SR), 
\item
[$\bm{r}_{\textbf{R}}(i)$,]
the direct ranks by the respondent (RR), 
\item
[$\bm{r}_{\bm{s}_{\textbf{R}}}(i)$,]
the ranks based on direct ratings by the respondent (RSR), 
\item
[$\bm{s}_{\textbf{PE}}(i)$,]
the ratings by the method of principal eigenvector (SPE), 
\item
[$\bm{r}_{\bm{s}_{\textbf{PE}}}(i)$,]
the ranks based on the principal eigenvector ratings (RSPE), 
\item
[$\bm{s}_{\textbf{GM}}(i)$,]
the ratings by the method of geometric mean (SGM), 
\item
[$\bm{r}_{\bm{s}_{\textbf{GM}}}(i)$,]
the ranks based on the geometric mean ratings (RSGM), 
\item
[$\bm{s}_{\textbf{CB}}(i)$,]
the best ratings by the method of log-Chebyshev approximation (SCB), 
\item
[$\bm{r}_{\bm{s}_{\textbf{CB}}}(i)$,]
the ranks based on the best ratings from the log-Chebyshev approximation (RSCB), 
\item
[$\bm{s}_{\textbf{CW}}(i)$,]
the worst ratings by the method of log-Chebyshev approximation (SCW),
\item
[$\bm{r}_{\bm{s}_{\textbf{CW}}}(i)$,]
the ranks based on the worst ratings from the log-Chebyshev approximation~(RSCW).
\end{description}

For respondents $i=1,2,3$ considered above as examples, we have the following results. For respondent 1, the vectors of ratings and ranks take the form
\begin{gather*}
\bm{r}_{\bm{s}_{\textbf{R}}}(1)
=
\left(
\begin{array}{c}
3
\\
1
\\
5
\\
4
\\
2
\\
6
\end{array}
\right),
\qquad
\bm{s}_{\textbf{PE}}(1)
=
\left(
\begin{array}{c}
0.3581
\\
0.6526
\\
0.1142
\\
0.1832
\\
1.0000
\\
0.0943
\end{array}
\right),
\qquad
\bm{r}_{\bm{s}_{\textbf{PE}}}(1)
=
\left(
\begin{array}{c}
3
\\
2
\\
5
\\
4
\\
1
\\
6
\end{array}
\right),
\\
\bm{s}_{\textbf{GM}}(1)
=
\left(
\begin{array}{c}
0.3588
\\
0.6680
\\
0.1190
\\
0.1906
\\
1.0000
\\
0.0981
\end{array}
\right),
\qquad
\bm{s}_{\textbf{CB}}(1)
=
\left(
\begin{array}{c}
0.3218
\\
0.6551
\\
0.1036
\\
0.1581
\\
1.0000	
\\
0.1018
\end{array}
\right),
\qquad
\bm{s}_{\textbf{CW}}(1)
=
\left(
\begin{array}{c}
0.3218
\\
0.6551
\\
0.1036
\\
0.1581
\\
1.0000
\\
0.1018
\end{array}
\right),
\\
\bm{r}_{\bm{s}_{\textbf{GM}}}(1)
=
\left(
\begin{array}{c}
3
\\
2
\\
5
\\
4
\\
1
\\
6
\end{array}
\right),
\qquad
\bm{r}_{\bm{s}_{\textbf{CB}}}(1)
=
\left(
\begin{array}{c}
3
\\
2
\\
5
\\
4
\\
1
\\
6
\end{array}
\right),
\qquad
\bm{r}_{\bm{s}_{\textbf{CW}}}(1)
=
\left(
\begin{array}{c}
3
\\
2
\\
5
\\
4
\\
1
\\
6
\end{array}
\right).
\end{gather*}

From the results obtained, it follows that the ranks of criteria provided by all methods for respondent 1 coincide and arrange the criteria in the following order of preference with the symbol $\succ$ used to indicate the preference relation:
\begin{equation*}
\mathcal{C}_{5}
\succ
\mathcal{C}_{2}
\succ
\mathcal{C}_{1}
\succ
\mathcal{C}_{4}
\succ
\mathcal{C}_{6}.
\end{equation*}

Note that the best and worst differentiating solutions lead the same vector of ratings, and thus the method of log-Chebyshev approximation has a unique solution in this case.

For respondent 2, the calculation results in the following vectors of ratings and ranks:
\begin{gather*}
\bm{r}_{\bm{s}_{\textbf{R}}}(2)
=
\left(
\begin{array}{c}
1
\\
2
\\
3
\\
5
\\
4
\\
6
\end{array}
\right),
\qquad
\bm{s}_{\textbf{PE}}(2)
=
\left(
\begin{array}{c}
0.8485
\\
1.0000
\\
0.4462
\\
0.3170
\\
0.3467
\\
0.1392
\end{array}
\right),
\qquad
\bm{r}_{\bm{s}_{\textbf{PE}}}(2)
=
\left(
\begin{array}{c}
2
\\
1
\\
3
\\
5
\\
4
\\
6
\end{array}
\right),
\\
\bm{s}_{\textbf{GM}}(2)
=
\left(
\begin{array}{c}
0.8754
\\
1.0000
\\
0.4292
\\
0.3184
\\
0.3645
\\
0.1434
\end{array}
\right),
\qquad
\bm{s}_{\textbf{CB}}(2)
=
\left(
\begin{array}{c}
0.7500
\\
1.0000
\\
0.4543
\\
0.2752
\\
0.2500
\\
0.1101
\end{array}
\right),
\qquad
\bm{s}_{\textbf{CW}}(2)
=
\left(
\begin{array}{c}
1.0000
\\
1.0000
\\
0.4543
\\
0.2752
\\
0.4543
\\
0.1667
\end{array}
\right),
\\
\bm{r}_{\bm{s}_{\textbf{GM}}}(2)
=
\left(
\begin{array}{c}
2
\\
1
\\
3
\\
5
\\
4
\\
6
\end{array}
\right),
\qquad
\bm{r}_{\bm{s}_{\textbf{CB}}}(2)
=
\left(
\begin{array}{c}
2
\\
1
\\
3
\\
4
\\
5
\\
6
\end{array}
\right),
\qquad
\bm{r}_{\bm{s}_{\textbf{CW}}}(2)
=
\left(
\begin{array}{c}
1
\\
2
\\
3
\\
5
\\
4
\\
6
\end{array}
\right).
\end{gather*}

We observe that the ranks of the criteria obtained by both methods of the principal vector and of geometric mean are the same and determine the order
\begin{equation*}
\mathcal{C}_{2}
\succ
\mathcal{C}_{1}
\succ
\mathcal{C}_{3}
\succ
\mathcal{C}_{5}
\succ
\mathcal{C}_{4}
\succ
\mathcal{C}_{6}.
\end{equation*}

Furthermore, we may consider that the ratings, which are provided by the worst differentiating vector of log-Chebyshev approximation, lead to the above ranking as well. Indeed, this method gives equal ratings to the first and second criteria to allow these criteria to be ranked in any order, one of which is the same as above.

The best and worst differentiating vectors of ratings obtained as a result of the log-Chebyshev approximation differ in some elements and can be coupled as a vector with interval values in the form
\begin{equation*}
\bm{s}_{\textbf{C}}(2)
=
\left(
\begin{array}{c}
0.7500\ \dots\ 1.0000 
\\
1.0000
\\
0.4543
\\
0.2752
\\
0.2500\ \dots\ 0.4543
\\
0.1101\ \dots\ 0.1667
\end{array}
\right).
\end{equation*}

The order provided by these vectors can be represented in combined form as 
\begin{equation*}
\mathcal{C}_{2}
\succeq
\mathcal{C}_{1}
\succ
\mathcal{C}_{3}
\succeq
\mathcal{C}_{5}
\parallel
\mathcal{C}_{4}
\succ
\mathcal{C}_{6},
\end{equation*}
where the symbol $\succeq$ denotes the weak preference relation (preferred or indifferent) and $\parallel$ indicates an undetermined preference.

Finally, the results obtained for respondent 3 are given by the vectors 
\begin{gather*}
\bm{r}_{\bm{s}_{\textbf{R}}}(3)
=
\left(
\begin{array}{c}
2
\\
1
\\
6
\\
4
\\
3
\\
5
\end{array}
\right),
\qquad
\bm{r}_{\bm{s}_{\textbf{PE}}}(3)
=
\left(
\begin{array}{c}
0.3773
\\
1.0000
\\
0.0930
\\
0.1630
\\
0.6213
\\
0.1631
\end{array}
\right),
\qquad
\bm{r}_{\bm{s}_{\textbf{PE}}}(3)
=
\left(
\begin{array}{c}
3
\\
1
\\
6
\\
5
\\
2
\\
4
\end{array}
\right),
\\
\bm{s}_{\textbf{GM}}(3)
=
\left(
\begin{array}{c}
0.3865
\\
1.0000
\\
0.0916
\\
0.1710
\\
0.6368
\\
0.1728
\end{array}
\right),
\qquad
\bm{s}_{\textbf{CB}}(3)
=
\left(
\begin{array}{c}
0.3798
\\
1.0000
\\
0.1082
\\
0.1755
\\
0.6163
\\
0.1755
\end{array}
\right),
\qquad
\bm{s}_{\textbf{CW}}(3)
=
\left(
\begin{array}{c}
0.3798
\\
1.0000
\\
0.1082
\\
0.1755
\\
0.6163
\\
0.2279
\end{array}
\right),
\\
\bm{r}_{\bm{s}_{\textbf{GM}}}(3)
=
\left(
\begin{array}{c}
3
\\
1
\\
6
\\
5
\\
2
\\
4
\end{array}
\right),
\qquad
\bm{r}_{\bm{s}_{\textbf{CB}}}(3)
=
\left(
\begin{array}{c}
3
\\
1
\\
6
\\
4
\\
2
\\
5
\end{array}
\right),
\qquad
\bm{r}_{\bm{s}_{\textbf{CW}}}(3)
=
\left(
\begin{array}{c}
3
\\
1
\\
6
\\
5
\\
2
\\
4
\end{array}
\right).
\end{gather*}

In this case, the solutions provided by the method of principal eigenvector and method of geometric mean as well as the worst differentiating solution by the log-Chebyshev approximation define the same ranks of criteria, which yields the order
\begin{equation*}
\mathcal{C}_{2}
\succ
\mathcal{C}_{5}
\succ
\mathcal{C}_{1}
\succ
\mathcal{C}_{6}
\succ
\mathcal{C}_{4}
\succ
\mathcal{C}_{3}.
\end{equation*}

The best and worst differentiating vectors of ratings differ only by the last element and can be combined as follows:
\begin{equation*}
\bm{s}_{\textbf{C}}(3)
=
\left(
\begin{array}{c}
0.3798
\\
1.0000
\\
0.1082
\\
0.1755
\\
0.6163
\\
0.1755\ \dots\ 0.2279
\end{array}
\right).
\end{equation*}

The order of criteria, which corresponds to these vectors, can be represented as
\begin{equation*}
\mathcal{C}_{2}
\succ
\mathcal{C}_{5}
\succ
\mathcal{C}_{1}
\succ
\mathcal{C}_{6}
\succeq
\mathcal{C}_{4}
\succ
\mathcal{C}_{3}.
\end{equation*}

It follows from the above examples that the results of determining the ranks of criteria, which are obtained from their pairwise comparisons by applying different methods may vary. In the next sections, we examine the difference between these results more closely by using appropriate techniques of data analysis.

\section{Statistical Analysis of Ratings and Ranks}
\label{S-SARR}

In this section, we compare ratings and ranks of criteria, obtained by different methods for all respondents participating in the survey. Ratings and ranks that are directly specified by respondents are compared with those which are based on three methods of evaluating criteria from pairwise comparisons. We examine the results of the principal eigenvector method, the geometric mean method, and the method of log-Chebyshev approximation. For the latter method, we consider two solution vectors, which best and worst differentiate the criteria if the solution is nonunique (both vectors coincide when the solution is unique).

To represent the obtained results in table form below, we use the following symbols:
\begin{description}[noitemsep, left=\parindent]
\item[SR,]
the direct ratings (scores) given by respondent (Scores by Respondent),
\item[RR,]
the direct ranks by respondent (Ranks by Respondent),
\item[RSR,]
the ranks from the direct ratings by respondent (Ranks from Scores by Respondent),
\item[SPE,]
the ratings by the method of principal eigenvector (Scores by Principal Eigenvector),
\item[RSPE,]
the ranks from the principal eigenvector ratings (Ranks from Scores by Principal~Eigenvector),
\item[SGM,]
the ratings by the method of geometric mean (Scores by Geometric Means),
\item[RSGM,]
the ranks from the geometric mean ratings (Ranks from Scores by Geometric~Means),
\item[SCB,]
the best differentiating ratings by the method of log-Chebyshev approximation (Scores by log-Chebyshev, Best),
\item[RSCB,]
the ranks from the best ratings from the log-Chebyshev approximation (Ranks from Scores by log-Chebyshev, Best),
\item[SCW,]
the worst differentiating ratings from the log-Chebyshev approximation (Scores by log-Chebyshev, Worst),
\item[RSCW,]
the ranks from the worst ratings from the log-Chebyshev approximation (Ranks from Scores by log-Chebyshev, Worst).
\end{description}

We start by counting the matches in the rank vectors obtained directly and determined through pairwise comparisons. The numbers of matches between rank vectors for each pair of the ranking methods are presented in Table~\ref{T-NMRV}.
\begin{specialtable}[ht]
\small
\caption{Number of matches of rank vectors.\label{T-NMRV}}
\begin{tabular*}{\hsize}{@{}@{\extracolsep{\fill}}ccccccc@{}}
\toprule
\textbf{Ranking} & \multicolumn{6}{c}{\textbf{Ranking Method}} \\
\cmidrule{2-7}
\textbf{Method} & \textbf{RR} & \textbf{RSR} & \textbf{RSPE} & \textbf{RSGM} & \textbf{RSCB} & \textbf{RSCW} \\
\midrule
RR & 202 & 28 & 56 & 56 & 59 & 56 \\
RSR & 28 & 202 & 20 & 20 & 18 & 23 \\
RSPE & 56 & 20 & 202 & 184 & 124 & 123 \\
RSGM & 56 & 20 & 184 & 202 & 124 & 125 \\
RSCB & 59 & 18 & 124 & 124 & 202 & 130 \\
RSCW & 56 & 23 & 123 & 125 & 130 & 202 \\
\bottomrule
\end{tabular*}
\end{specialtable}

The results given by Table~\ref{T-NMRV} show a low level of matches (encountered for 28 out of 202 respondents) between the vectors of direct ranks (RR) and ranks determined from ratings that are specified by respondents (RSR). We can explain this inconsistency by greater variability in determining ratings by respondents, which allows one to assign the same ratings to several criteria. As a result, the ranks obtained from ratings may differ from the direct ranking, when the respondent specifies an individual rank for each criterion to indicate the number (position) of the criterion in the ranking order.

The degree of consistency of the direct ranking (RR) with the ranking on the basis of ratings obtained from pairwise comparisons (RSPE, RSGM, RSCB, RSCW) is better given by 56--59 matches of rank vectors, which are more than a quarter of respondents. Note that the maximum number of matches, which is equal to 59, is provided by the best differentiating solution of log-Chebyshev approximation (RSCB).

We note that the results of ranking criteria from pairwise comparisons demonstrate a higher level of consistency between the methods under consideration. Specifically, the largest number of matches encountered for 184 (more than 90\% of) respondents is given by the outcome from the methods of principal eigenvector (RSPE) and of geometric mean (RSGM). The number of matches between the rank vectors produced by these methods and the vectors obtained using log-Chebyshev approximation (RSCB, RSCW) is within the range 123--125 (more than 60\% of respondents). Finally, the rank vectors, which are obtained from the best and worst differentiating solutions found by log-Chebyshev approximation, coincide for 130 respondents.

Below, to analyze and compare the results more thoroughly, we use appropriate statistical techniques of parameter estimation and correlation analysis. As preliminary analysis of data, we investigate if parameters of respondents, such as age and sex, may affect the results of the evaluation of criteria. Furthermore, we compare the results of evaluating ratings and ranks obtained from pairwise comparisons by different methods with each other and with the ratings and ranks directly provided by respondents.

\subsection{Preliminary Analysis of Data}

In the analysis of indirect methods that produce results from pairwise comparison data, it seems reasonable to consider the ratings and ranks of criteria, which are directly specified by respondents, as some basis for comparison. Since these ratings and ranks are known for all respondents, we can try to select those respondents for whom the direct ranks and ranks given by direct ratings appear to be closest to each other. Considering that a part of respondents may demonstrate a higher level of consistency in the direct judgment, one can expect that these respondents are able to assess pairwise comparison of criteria more adequately and definitely.

To examine how the consistency between direct estimates of ratings and ranks may affect the results based on pairwise comparisons, we consider groups of respondents with different degrees of consistency. As a measure of inconsistency (difference) between rank vectors, we apply the Chebyshev metric (the maximum absolute componentwise difference), which is calculated for two $n$-vectors $\bm{a}=(a_{j})$ and $\bm{b}=(b_{j})$ as
\begin{equation*}
d(\bm{a},\bm{b})
=
\max_{1\leq j\leq n}|a_{j}-b_{j}|.
\end{equation*}

We consider a series of nested subsets (groups) $R_{0}\subset R_{1}\subset\cdots\subset R_{5}$, where $R_{i}$ denotes the set of respondents for whom the Chebyshev metric between the vectors of direct ranks and ranks obtained from the direct ratings is less or equal to $i=0,1,\ldots,5$. Table~\ref{T-GRADBRV} shows how the total number of respondents in groups changes with increase in the upper bound for Chebyshev metric, which indicates the maximum degree of difference in groups. For each group, the percentage of respondents who visited the country and male/female percentage are also included.
\begin{specialtable}[ht]
\small
\caption{Groups of respondents according to difference between rank vectors.\label{T-GRADBRV}}
\begin{tabular*}{\hsize}{@{}@{\extracolsep{\fill}}ccccc@{}}
\toprule
								& \textbf{Maximum} 				& \textbf{Number of} 		& \textbf{Percentage of} 				& \textbf{Male/} \\
\textbf{Group}	& \textbf{Component-Wise} & \textbf{Respondents}	& \textbf{Respondents Who}	 		& \textbf{Female} \\
								& \textbf{Difference}			& \textbf{in the Group}	& \textbf{Visited the Country}	& \textbf{Percentage} \\
\midrule		
$R_{0}$					& $0$ 		 						& 28 									& 35 				 									& 33/67 \\
$R_{1}$					& $1$			 								& 95 									& 19 				 									& 32/68 \\
$R_{2}$					&	$2$ 				 						& 167 									& 18 			 									& 30/70 \\
$R_{3}$					& $3$ 		 								& 194 									& 17 				 									& 28/72 \\
$R_{4}$					& $4$ 		 								& 201 									& 18 				 									& 28/72 \\
$R_{5}$					& $5$ 		 								& 202 									& 19 				 									& 28/72 \\
\bottomrule
\end{tabular*}
\end{specialtable}

It follows from the data presented in the table that the direct ranks and the ranks obtained from direct ratings completely coincide for 28 respondents from the group $R_{0}$, who may be considered as the most accurate evaluators. However, these respondents form a rather small part of all respondents involved in the survey, which makes it unreasonable to take this group as a good representative of the entire sample.

Furthermore, we note that the respondents in the groups $R_{1},\ldots,R_{5}$ with less consistent direct judgments demonstrate almost the same percentage of those who visited the country (17--19\%). This apparently indicates that for the most respondents, there is no relationship between visiting the country and accuracy of judgments. At the same time, for the group $R_{0}$, this percentage increases to 35\%, which is in line with the idea that the accuracy of judgments should increase as respondents gain experience of staying in the country. A~conclusion about the lack of systematic dependence of the accuracy of judgments on the gender and age of respondents can also be drawn for all groups.

To complete the analysis of the possible influence of the accuracy of direct judgments about ratings and ranks on the other survey results, we estimate correlations between the rank vectors directly provided by respondents and the vectors derived from pairwise comparisons. To measure correlation, we use the Kendall rank correlation coefficient, which is given for two $n$-vectors $\bm{a}=(a_{j})$ and $\bm{b}=(b_{j})$ by the formula
\begin{equation*}
\tau(\bm{a},\bm{b})
=
\frac{2}{n(n-1)}
\sum_{1\leq i<j\leq n}
\mathop\mathrm{sgn}(a_{i}-a_{j})\mathop\mathrm{sgn}(b_{i}-b_{j}).
\end{equation*}

For each group $R_{0}$, $R_{1}$ and $R_{5}$, we concatenate the vectors of ratings obtained for all respondents in the group to form a single vector of ratings for each method. The estimates of the correlation coefficients between concatenated vectors for the groups $R_{0}$, $R_{1}$ and $R_{5}$ are given in Tables~\ref{T-CGR0}--\ref{T-CGR5} (the estimates for groups $R_{2}$, $R_{3}$ and $R_{4}$ do not differ significantly from the results for $R_{1}$ and $R_{5}$ and hence are omitted). 
\begin{specialtable}[ht]
\small
\caption{Correlation of rank vectors for the group $R_{0}$.\label{T-CGR0}}
\begin{tabular*}{\hsize}{@{}@{\extracolsep{\fill}}cccccc@{}}
\toprule
\textbf{Ranking} & \multicolumn{5}{c}{\textbf{Ranking Method}} \\
\cmidrule{2-6}
\textbf{Method} & \textbf{RR} & \textbf{RSPE} & \textbf{RSGM} & \textbf{RSCB} & \textbf{RSCW} \\
\midrule
RR & 1.0000 & 0.8482 & 0.8482 & 0.8560 & 0.8645 \\
RSPE & 0.8482 & 1.0000 & 1.0000 & 0.9558 & 0.9557 \\
RSGM & 0.8482 & 1.0000 & 1.0000 & 0.9557 & 0.9557 \\
RSCB & 0.8560 & 0.9558 & 0.9558 & 1.0000 & 0.9717 \\
RSCW & 0.8645 & 0.9557 & 0.9557 & 0.9717 & 1.0000 \\
\bottomrule
\end{tabular*}
\end{specialtable}

\begin{specialtable}[ht]
\small
\caption{Correlation of rank vectors for the group $R_{1}$.\label{T-CGR1}}
\begin{tabular*}{\hsize}{@{}@{\extracolsep{\fill}}cccccc@{}}
\toprule
\textbf{Ranking} & \multicolumn{5}{c}{\textbf{Ranking Method}} \\
\cmidrule{2-6}
\textbf{Method} & \textbf{RR} & \textbf{RSPE} & \textbf{RSGM} & \textbf{RSCB} & \textbf{RSCW} \\
\midrule
RR & 1.0000 & 0.8184 & 0.8234 & 0.8108 & 0.8221 \\
RSPE & 0.8184 & 1.0000 & 0.9943 & 0.9450 & 0.9183 \\
RSGM & 0.8234 & 0.9943 & 1.0000 & 0.9429 & 0.9223 \\
RSCB & 0.8108 & 0.9450 & 0.9429 & 1.0000 & 0.9234 \\
RSCW & 0.8221 & 0.9183 & 0.9223 & 0.9234 & 1.0000 \\
\bottomrule
\end{tabular*}
\end{specialtable}

\begin{specialtable}[ht]
\small
\caption{Correlation of rank vectors for the group $R_{5}$.\label{T-CGR5}}
\begin{tabular*}{\hsize}{@{}@{\extracolsep{\fill}}cccccc@{}}
\toprule
\textbf{Ranking} & \multicolumn{5}{c}{\textbf{Ranking Method}} \\
\cmidrule{2-6}
\textbf{Method} & \textbf{RR} & \textbf{RSPE} & \textbf{RSGM} & \textbf{RSCB} & \textbf{RSCW} \\
\midrule
RR & 1.0000 & 0.7939 & 0.7940 & 0.7876 & 0.7839 \\
RSPE & 0.7939 & 1.0000 & 0.9862 & 0.9312 & 0.9186 \\
RSGM & 0.7940 & 0.9862 & 1.0000 & 0.9290 & 0.9204 \\
RSCB & 0.7876 & 0.9312 & 0.9290 & 1.0000 & 0.9144 \\
RSCW & 0.7839 & 0.9186 & 0.9204 & 0.9144 & 1.0000 \\
\bottomrule
\end{tabular*}
\end{specialtable}

The data presented in these tables show that for all groups of respondents, the correlations vary within small limits and maintain the relative order of magnitude for different methods. This allows us to conclude that the change in accuracy of direct estimates does not significantly affect the results of ranking based on pairwise comparisons.

In the subsequent analysis, we take into account the above results and examine all respondents together rather than divide them into groups according to the accuracy of direct judgments or parameters of respondents.

\subsection{Comparison of Ratings}

We now examine results of evaluating ratings, obtained by different methods from pairwise comparisons and specified directly by respondents. First, we calculate the means and standard deviations of vectors of ratings for each rating method. The mean vectors can be considered as the most representative (typical) vectors over all respondents for each method to give some general idea of the absolute ratings of criteria in hotel selection. The standard deviations describing the spread of the vectors can characterize both the differences in the perception of the relative importance of criteria by respondents and the ability of the methods to reinforce or smooth these differences.

The results of calculating the mean vectors presented in Table~\ref{T-MVRC} show that on average, criterion $2$ is always rated higher than the others, next come criterion $1$ and then $5$. The other criteria have lower mean ratings, which give the order $(3,4,6)$ for all methods except the direct rating with the order $(6,3,4)$. Examination of the output of individual respondents indicates that the second criterion is evaluated above the others for more than 60\% respondents. Specifically, this result encounters in the direct rating (SR) for 165~respondents, in both methods of principal eigenvector (SPE) and geometric mean (SGM) for 125 respondents, and in the best and worst differentiating solutions of log-Chebyshev approximation (SCB and SCW) for 123 and 147 respondents respectively.
\begin{specialtable}[ht]
\small
\caption{Mean vectors of rating of criteria.\label{T-MVRC}}
\begin{tabular*}{\hsize}{@{}@{\extracolsep{\fill}}cccccc@{}}
\toprule
\textbf{Criterion} & \multicolumn{5}{c}{\textbf{Rating Method}} \\
\cmidrule{2-6}
\textbf{Number} & \textbf{SR} & \textbf{SPE} & \textbf{SGM} & \textbf{SCB} & \textbf{SCW} \\
\midrule
1 & 0.7480 & 0.6268 & 0.6272 & 0.5693 & 0.6313 \\
2 & 0.9428 & 0.8877 & 0.8915 & 0.8504 & 0.9025 \\
3 & 0.5784 & 0.3238 & 0.3231 & 0.2861 & 0.3511 \\
4 & 0.5588 & 0.3044 & 0.3053 & 0.2735 & 0.3291 \\
5 & 0.7331 & 0.5279 & 0.5303 & 0.4910 & 0.5581 \\
6 & 0.6149 & 0.2543 & 0.2558 & 0.2252 & 0.2879 \\
\bottomrule
\end{tabular*}
\end{specialtable}

It follows from Table~\ref{T-MVRC} that the direct ratings by respondents are sufficiently different from the other results. This can be explained by a rather primitive scale used in the survey, which includes only 5 points. At the same time, the evaluation of mean vectors for all indirect methods of rating criteria gives quite close results. We note that the mean vectors of the log-Chebyshev approximation appear as lower and upper boundaries for the mean vectors for both methods of principal eigenvector and geometric mean. 

To measure the overall variability of ratings for all methods, we use the average standard deviation over criteria and a total (vector) standard deviation defined as the square root of the sum of squares of standard deviations for each criterion. The results of calculation are given in Table~\ref{T-VEVR} and say, in particular, that the total standard deviations for the methods of principal eigenvector and geometric mean almost coincide, and both methods have slightly lower variability of ratings than the log-Chebyshev approximation.
\begin{specialtable}[ht]
\small
\caption{Variability estimates for vectors of ratings.\label{T-VEVR}}
\begin{tabular*}{\hsize}{@{}@{\extracolsep{\fill}}cccccc@{}}
\toprule
\textbf{Variability} & \multicolumn{5}{c}{\textbf{Rating Method}} \\
\cmidrule{2-6}
\textbf{Measure} & \textbf{SR} & \textbf{SPE} & \textbf{SGM} & \textbf{SCB} & \textbf{SCW} \\
\midrule
Average deviation & 0.2037 & 0.2202 & 0.2198 & 0.2252 & 0.2339 \\
Total deviation & 0.5046 & 0.5477 & 0.5478 & 0.5613 & 0.5796 \\
\bottomrule
\end{tabular*}
\end{specialtable}

Table~\ref{T-SDREC} demonstrates the standard deviation of ratings for each pair of criteria and methods. We note that criteria 2 and 6 most often have the smallest or the next smallest standard deviation for all methods, which suggests that the real assessment of these criteria differs the least from one respondent to another. 
\begin{specialtable}[ht]
\small
\caption{Standard deviation of ratings of each criterion.\label{T-SDREC}}
\begin{tabular*}{\hsize}{@{}@{\extracolsep{\fill}}cccccc@{}}
\toprule
\textbf{Criterion} & \multicolumn{5}{c}{\textbf{Rating Method}} \\
\cmidrule{2-6}
\textbf{Number} & \textbf{SR} & \textbf{SPE} & \textbf{SGM} & \textbf{SCB} & \textbf{SCW} \\
\midrule
1 &	0.2188 & 0.2578 & 0.2611 & 0.2685 & 0.2633 \\
2 & 0.1379 & 0.1837 & 0.1807 & 0.2175 & 0.1879 \\
3 & 0.2255 & 0.2300 & 0.2279 & 0.2163 & 0.2426 \\
4 & 0.2248 & 0.2052 & 0.2039 & 0.1963 & 0.2150 \\
5 & 0.1978 & 0.2769 & 0.2798 & 0.2889 & 0.2914 \\
6 & 0.2173 & 0.1676 & 0.1656 & 0.1639 & 0.2032 \\
\bottomrule
\end{tabular*}
\end{specialtable}

Next, we turn to the estimates of correlation between vectors of ratings produced by the direct and indirect assessment methods. For each method, we concatenate the vectors of all respondents into a single vector of ratings, and then estimate the Pearson correlation coefficients between the concatenated vectors. The correlation coefficients calculated for all pairs of methods are given in Table~\ref{T-CBVRAPM}.
\begin{specialtable}[ht]
\small
\caption{Correlation of vectors of ratings for all pairs of methods.\label{T-CBVRAPM}}
\begin{tabular*}{\hsize}{@{}@{\extracolsep{\fill}}cccccc@{}}
\toprule
\textbf{Rating} & \multicolumn{5}{c}{\textbf{Rating Method}} \\
\cmidrule{2-6}
\textbf{Method} & \textbf{SR} & \textbf{SPE} & \textbf{SGM} & \textbf{SCB} & \textbf{SCW} \\
\midrule
SR	& 1.0000 & 0.6942 & 0.6918 & 0.6619 & 0.6889 \\
SPE & 0.6942 & 1.0000 & 0.9973 & 0.9613 & 0.9622 \\
SGM & 0.6918 & 0.9973 & 1.0000 & 0.9631 & 0.9576 \\
SCB & 0.6619 & 0.9613 & 0.9631 & 1.0000 & 0.9164 \\
SCW & 0.6889 & 0.9622 & 0.9576 & 0.9164 & 1.0000 \\
\bottomrule
\end{tabular*}
\end{specialtable}

The results presented in this table show that the correlation between ratings obtained from pairwise comparisons is close to one for all pairs of methods, which means that any two concatenated vectors of ratings are close to be linearly dependent (collinear). Specifically, the vectors found by the methods of principal eigenvector and geometric mean correlate at the level of $0.9973$, and thus can be considered poorly distinguishable from the correlation point of view.

We now consider the correlations of ratings produced by all methods for each individual criterion. Tables~\ref{T-CRC1APM}--\ref{T-CRC6APM} demonstrate the estimates for criteria 1, 2 and 6 (the results for criteria 3, 4 and 5 are similar to those for criterion 1 and thus omitted). 
\begin{specialtable}[ht]
\small
\caption{Correlation of ratings of criterion 1 for all pairs of methods.\label{T-CRC1APM}}
\begin{tabular*}{\hsize}{@{}@{\extracolsep{\fill}}cccccc@{}}
\toprule
\textbf{Rating} & \multicolumn{5}{c}{\textbf{Rating Method}} \\
\cmidrule{2-6}
\textbf{Method} & \textbf{SR} & \textbf{SPE} & \textbf{SGM} & \textbf{SCB} & \textbf{SCW} \\
\midrule
SR & 1.0000 & 0.5815 & 0.5639 & 0.5151 & 0.5758 \\
SPE & 0.5815 & 1.0000 & 0.9949 & 0.9215 & 0.9408 \\
SGM & 0.5639 & 0.9949 & 1.0000 & 0.9250 & 0.9220 \\
SCB & 0.5151 & 0.9215 & 0.9250 & 1.0000 & 0.8723 \\
SCW & 0.5758 & 0.9408 & 0.9220 & 0.8723 & 1.0000 \\
\bottomrule
\end{tabular*}
\end{specialtable}

The results obtained for each criterion indicate that there is a high correlation at the level of $0.9929\dots0.9959$ between the ratings produced by the method of principal eigenvector and the method of geometric mean. 
\begin{specialtable}[ht]
\small
\caption{Correlation of ratings of criterion 2 for all pairs of methods.\label{T-CRC2APM}}
\begin{tabular*}{\hsize}{@{}@{\extracolsep{\fill}}cccccc@{}}
\toprule
\textbf{Rating} & \multicolumn{5}{c}{\textbf{Rating Method}} \\
\cmidrule{2-6}
\textbf{Method} & \textbf{SR} & \textbf{SPE} & \textbf{SGM} & \textbf{SCB} & \textbf{SCW} \\
\midrule
SR & 1.0000 & 0.4243 & 0.4314 & 0.3778 & 0.4348 \\
SPE & 0.4243 & 1.0000 & 0.9959 & 0.8887 & 0.9217 \\
SGM & 0.4314 & 0.9959 & 1.0000 & 0.8899 & 0.9239 \\
SCB & 0.3778 & 0.8887 & 0.8899 & 1.0000 & 0.8004 \\
SCW & 0.4348 & 0.9217 & 0.9239 & 0.8004 & 1.0000 \\
\bottomrule
\end{tabular*}
\end{specialtable}

Furthermore, the ratings provided by the best differentiating solutions of log-Chebyshev approximation are usually more correlated with the ratings by the method of geometric mean ($0.8899\dots0.9543$), whereas the ratings from the worst differentiating solutions are more correlated with those from the method of principal eigenvector ($0.8938\dots0.9408$). Finally, the correlation between ratings from the best and worst differentiating solutions are somewhat lower and lay in the range $0.7776\dots0.8752$.
\begin{specialtable}[ht]
\small
\caption{Correlation of ratings of criterion 6 for all pairs of methods.\label{T-CRC6APM}}
\begin{tabular*}{\hsize}{@{}@{\extracolsep{\fill}}cccccc@{}}
\toprule
\textbf{Rating} & \multicolumn{5}{c}{\textbf{Rating Method}} \\
\cmidrule{2-6}
\textbf{Method} & \textbf{SR} & \textbf{SPE} & \textbf{SGM} & \textbf{SCB} & \textbf{SCW} \\
\midrule
SR & 1.0000 & 0.4926 & 0.4873 & 0.4569	& 0.4842 \\
SPE & 0.4926 & 1.0000 & 0.9929 & 0.9301 & 0.8938 \\
SGM & 0.4873 & 0.9929 & 1.0000 & 0.9189 & 0.8911 \\
SCB & 0.4569 & 0.9301 & 0.9189 & 1.0000 & 0.7776 \\
SCW & 0.4842 & 0.8938 & 0.8911 & 0.7776 & 1.0000 \\
\bottomrule
\end{tabular*}
\end{specialtable}

To summarize the results of the comparison of ratings, we can conclude that for the survey data under study, all indirect methods of rating criteria consistently produce highly correlated vectors of ratings, with the highest correlation very close to one between the results of the methods of principal eigenvector and geometric mean. The good agreement between the results of indirect methods can be explained by the natural origin of the initial data obtained as results of human judgment, which provide rationale for more consistent pairwise comparisons than artificial simulated data.

In conclusion, we note low correlations between direct and indirect ratings, which are within the range $0.3778\dots0.5993$. This relative disagreement may serve as an illustration of the known fact that it may be difficult for typical respondents to make a correct judgment when more than two alternatives are simultaneously evaluated, whereas they normally assess pairwise comparisons more accurately.

\subsection{Comparison of Ranks}

We start with the ranks derived from ratings given by the mean vectors in Table~\ref{T-MVRC}. The mean vector of ratings that are directly provided by respondents yields the rank vector $(2,1,5,6,3,4)$, which arranges the criteria in the following order: 
\begin{equation*}
\mathcal{C}_{2}\succ\mathcal{C}_{1}\succ\mathcal{C}_{5}\succ\mathcal{C}_{6}\succ\mathcal{C}_{3}\succ\mathcal{C}_{4}.
\end{equation*}

To see how frequently this order coincides with the order derived from the vectors of ratings of individual respondents, we calculate the Hamming distance between (the number of different elements in) the corresponding rank vectors. The results of counting the number of rank vectors, which are obtained from vectors of ratings for all methods, within fixed distances from the vector $(2,1,5,6,3,4)$ are given in Table~\ref{T-NRVFD215634}.
\begin{specialtable}[ht]
\small
\caption{Number of vectors within fixed distance from $(2,1,5,6,3,4)$.\label{T-NRVFD215634}}
\begin{tabular*}{\hsize}{@{}@{\extracolsep{\fill}}cccccc@{}}
\toprule
\textbf{Maximum} & \multicolumn{5}{c}{\textbf{Rating Method}} \\
\cmidrule{2-6}
\textbf{Hamming Distance} & \textbf{RSR} & \textbf{RSPE} & \textbf{RSGM} & \textbf{RSCB} & \textbf{RSCW} \\
\midrule
0 & 5 & 3 & 3 & 4 & 4 \\
1 & 5 & 3 & 3 & 4 & 4 \\
2 & 23 & 35 & 32 & 35 & 31 \\
3 & 61 & 68 & 65 & 66 & 64 \\
4 & 122 & 135 & 130 & 123 & 132 \\
5 & 178 & 185 & 186 & 186 & 188 \\
6 & 202 & 202 & 202 & 202 & 202 \\
\bottomrule
\end{tabular*}
\end{specialtable}

The mean vectors of ratings derived from pairwise comparisons by different methods imply a common rank vector $(2,1,5,3,4,6)$, which corresponds to the order
\begin{equation*}
\mathcal{C}_{2}\succ\mathcal{C}_{1}\succ\mathcal{C}_{5}\succ\mathcal{C}_{3}\succ\mathcal{C}_{4}\succ\mathcal{C}_{6}.
\end{equation*}

The numbers of rank vectors within fixed distances from the vector $(2,1,5,3,4,6)$ are shown in Table~\ref{T-NRVFD215346}.
\begin{specialtable}[ht]
\small
\caption{Number of vectors within fixed distance from $(2,1,5,3,4,6)$.\label{T-NRVFD215346}}
\begin{tabular*}{\hsize}{@{}@{\extracolsep{\fill}}cccccc@{}}
\toprule
\textbf{Maximum} & \multicolumn{5}{c}{\textbf{Rating Method}} \\
\cmidrule{2-6}
\textbf{Hamming Distance} & \textbf{RSR} & \textbf{RSPE} & \textbf{RSGM} & \textbf{RSCB} & \textbf{RSCW} \\
\midrule
0 & 6 & 4 & 4 & 7 & 5 \\
1 & 6 & 4 & 4 & 7 & 5 \\
2 & 37 & 40 & 38 & 48 & 41 \\
3 & 66 & 64 & 63 & 69 & 63 \\
4 & 121 & 133 & 132 & 130 & 129 \\
5 & 175 & 180 & 179 & 183 & 184 \\
6 & 202 & 202 & 202 & 202 & 202 \\
\bottomrule
\end{tabular*}
\end{specialtable}

Both Tables~\ref{T-NRVFD215634} and \ref{T-NRVFD215346} demonstrate low correspondence between the rank vectors derived from the mean vectors of ratings and the individual rank vectors for each respondent. In this case, we cannot consider the ranks based on the mean vectors as well justified, and need further examination of the individual rank vectors to find those vectors that occur more frequently for each method. The results of determining the five most common rank vectors for each method are presented in Table~\ref{T-MORV}.
\begin{specialtable}[ht]
\small
\caption{Most frequently occurred rank vectors.\label{T-MORV}}
\begin{tabular*}{\hsize}{@{}@{\extracolsep{\fill}}cccc@{}}
\toprule
\textbf{Ranking Method} & \textbf{Rank Vector} & \textbf{Order of Criteria} & \textbf{Number of Occurrences} \\
\midrule
 & $(1,2,5,3,4,6)$ & $\mathcal{C}_{1}\succ\mathcal{C}_{2}\succ\mathcal{C}_{5}\succ\mathcal{C}_{3}\succ\mathcal{C}_{4}\succ\mathcal{C}_{6}$ & 10 \\
 & $(2,1,5,3,6,4)$ & $\mathcal{C}_{2}\succ\mathcal{C}_{1}\succ\mathcal{C}_{5}\succ\mathcal{C}_{3}\succ\mathcal{C}_{6}\succ\mathcal{C}_{4}$ & 9 \\
RSR & $(1,2,3,5,6,4)$ & $\mathcal{C}_{1}\succ\mathcal{C}_{2}\succ\mathcal{C}_{3}\succ\mathcal{C}_{5}\succ\mathcal{C}_{6}\succ\mathcal{C}_{4}$ & 8 \\
 & $(2,5,6,1,3,4)$ & $\mathcal{C}_{2}\succ\mathcal{C}_{5}\succ\mathcal{C}_{6}\succ\mathcal{C}_{1}\succ\mathcal{C}_{3}\succ\mathcal{C}_{4}$ & 7 \\
 & $(1,2,3,4,5,6)$ & $\mathcal{C}_{1}\succ\mathcal{C}_{2}\succ\mathcal{C}_{3}\succ\mathcal{C}_{4}\succ\mathcal{C}_{5}\succ\mathcal{C}_{6}$ & 6 \\
\midrule
 & $(2,1,5,3,6,4)$ & $\mathcal{C}_{2}\succ\mathcal{C}_{1}\succ\mathcal{C}_{5}\succ\mathcal{C}_{3}\succ\mathcal{C}_{6}\succ\mathcal{C}_{4}$ & 12 \\
 & $(2,1,5,4,6,3)$ & $\mathcal{C}_{2}\succ\mathcal{C}_{1}\succ\mathcal{C}_{5}\succ\mathcal{C}_{4}\succ\mathcal{C}_{6}\succ\mathcal{C}_{3}$ & 8 \\
RSPE & $(2,1,5,4,3,6)$ & $\mathcal{C}_{2}\succ\mathcal{C}_{1}\succ\mathcal{C}_{5}\succ\mathcal{C}_{4}\succ\mathcal{C}_{3}\succ\mathcal{C}_{6}$ & 7 \\
 & $(2,5,1,3,6,4)$ & $\mathcal{C}_{2}\succ\mathcal{C}_{5}\succ\mathcal{C}_{1}\succ\mathcal{C}_{3}\succ\mathcal{C}_{6}\succ\mathcal{C}_{4}$ & 6 \\
 & $(2,1,4,5,6,3)$ & $\mathcal{C}_{2}\succ\mathcal{C}_{1}\succ\mathcal{C}_{4}\succ\mathcal{C}_{5}\succ\mathcal{C}_{6}\succ\mathcal{C}_{3}$ & 6 \\
\midrule
 & $(2,1,5,3,6,4)$ & $\mathcal{C}_{2}\succ\mathcal{C}_{1}\succ\mathcal{C}_{5}\succ\mathcal{C}_{3}\succ\mathcal{C}_{6}\succ\mathcal{C}_{4}$ & 12 \\
 & $(2,1,5,4,6,3)$ & $\mathcal{C}_{2}\succ\mathcal{C}_{1}\succ\mathcal{C}_{5}\succ\mathcal{C}_{4}\succ\mathcal{C}_{6}\succ\mathcal{C}_{3}$ & 8 \\
RSGM & $(2,1,4,5,6,3)$ & $\mathcal{C}_{2}\succ\mathcal{C}_{1}\succ\mathcal{C}_{4}\succ\mathcal{C}_{5}\succ\mathcal{C}_{6}\succ\mathcal{C}_{3}$ & 7 \\
 & $(1,2,5,4,3,6)$ & $\mathcal{C}_{1}\succ\mathcal{C}_{2}\succ\mathcal{C}_{5}\succ\mathcal{C}_{4}\succ\mathcal{C}_{3}\succ\mathcal{C}_{6}$ & 7 \\
 & $(2,1,3,5,6,4)$ & $\mathcal{C}_{2}\succ\mathcal{C}_{1}\succ\mathcal{C}_{3}\succ\mathcal{C}_{5}\succ\mathcal{C}_{6}\succ\mathcal{C}_{4}$ & 6 \\
\midrule
 & $(2,1,5,3,6,4)$ & $\mathcal{C}_{2}\succ\mathcal{C}_{1}\succ\mathcal{C}_{5}\succ\mathcal{C}_{3}\succ\mathcal{C}_{6}\succ\mathcal{C}_{4}$ & 11 \\
 & $(2,5,1,3,4,6)$ & $\mathcal{C}_{2}\succ\mathcal{C}_{5}\succ\mathcal{C}_{1}\succ\mathcal{C}_{3}\succ\mathcal{C}_{4}\succ\mathcal{C}_{6}$ & 7 \\
RSCB & $(2,1,5,3,4,6)$ & $\mathcal{C}_{2}\succ\mathcal{C}_{1}\succ\mathcal{C}_{5}\succ\mathcal{C}_{3}\succ\mathcal{C}_{4}\succ\mathcal{C}_{6}$ & 7 \\
 & $(2,1,5,4,3,6)$ & $\mathcal{C}_{2}\succ\mathcal{C}_{1}\succ\mathcal{C}_{5}\succ\mathcal{C}_{4}\succ\mathcal{C}_{3}\succ\mathcal{C}_{6}$ & 7 \\
 & $(1,2,5,4,3,6)$ & $\mathcal{C}_{1}\succ\mathcal{C}_{2}\succ\mathcal{C}_{5}\succ\mathcal{C}_{4}\succ\mathcal{C}_{3}\succ\mathcal{C}_{6}$ & 6 \\
\midrule
 & $(2,1,5,3,6,4)$ & $\mathcal{C}_{2}\succ\mathcal{C}_{1}\succ\mathcal{C}_{5}\succ\mathcal{C}_{3}\succ\mathcal{C}_{6}\succ\mathcal{C}_{4}$ & 9 \\
 & $(2,5,1,3,4,6)$ & $\mathcal{C}_{2}\succ\mathcal{C}_{5}\succ\mathcal{C}_{1}\succ\mathcal{C}_{3}\succ\mathcal{C}_{4}\succ\mathcal{C}_{6}$ & 8 \\
RSCW & $(2,1,5,4,3,6)$ & $\mathcal{C}_{2}\succ\mathcal{C}_{1}\succ\mathcal{C}_{5}\succ\mathcal{C}_{4}\succ\mathcal{C}_{3}\succ\mathcal{C}_{6}$ & 8 \\
 & $(1,2,5,4,3,6)$ & $\mathcal{C}_{1}\succ\mathcal{C}_{2}\succ\mathcal{C}_{5}\succ\mathcal{C}_{4}\succ\mathcal{C}_{3}\succ\mathcal{C}_{6}$ & 8 \\
 & $(2,5,1,4,3,6)$ & $\mathcal{C}_{2}\succ\mathcal{C}_{5}\succ\mathcal{C}_{1}\succ\mathcal{C}_{4}\succ\mathcal{C}_{3}\succ\mathcal{C}_{6}$ & 6 \\
\bottomrule
\end{tabular*}
\end{specialtable}

It follows from Table~\ref{T-MORV} that the most frequently occurred rank vector derived from ratings for all indirect methods and the second most frequent for the direct ratings by respondents is $(2,1,5,3,6,4)$. This rank vector yields the order of criteria defined as
\begin{equation*}
\mathcal{C}_{2}\succ\mathcal{C}_{1}\succ\mathcal{C}_{5}\succ\mathcal{C}_{3}\succ\mathcal{C}_{6}\succ\mathcal{C}_{4},
\end{equation*}
and seems to provide a more accurate assessment of the ranks of criteria for the set of respondents participating in the survey.

We conclude with the results in Table~\ref{T-CRVAM}, which demonstrate the Kendall correlation between concatenated rank vectors obtained by direct and indirect methods (this table actually presents the same results as in Table~\ref{T-CGR5} used in the preliminary analysis of data). 
\begin{specialtable}[ht]
\small
\caption{Correlations of rank vectors for all methods.\label{T-CRVAM}}
\begin{tabular*}{\hsize}{@{}@{\extracolsep{\fill}}cccccc@{}}
\toprule
\textbf{Ranking} & \multicolumn{5}{c}{\textbf{Ranking Method}} \\
\cmidrule{2-6}
\textbf{Method} & \textbf{RR} & \textbf{RSPE} & \textbf{RSGM} & \textbf{RSCB} & \textbf{RSCW} \\
\midrule
RR & 1.0000 & 0.7939 & 0.7940 & 0.7876 & 0.7839 \\
RSPE & 0.7939 & 1.0000 & 0.9862 & 0.9312 & 0.9186 \\
RSGM & 0.7940 & 0.9862 & 1.0000 & 0.9290 & 0.9204 \\
RSCB & 0.7876 & 0.9312 & 0.9290 & 1.0000 & 0.9144 \\
RSCW & 0.7839 & 0.9186 & 0.9204 & 0.9144 & 1.0000 \\
\bottomrule
\end{tabular*}
\end{specialtable}

As for the correlation of vectors of ratings obtained from pairwise comparisons, the corresponding rank vectors are highly correlated for all methods. As before, the methods of principal eigenvector and geometric mean have the highest correlation of $0.9863$, which is very close to one. At the same time, the correlation of rank vectors obtained by the log-Chebyshev approximation and by both methods of principal eigenvector and geometric mean increase to $0.9187\dots0.9313$. This shows that in terms of rank vectors, the difference between the results of all methods remains quite small, and therefore, the results of one method can be used to validate the results of the other methods.

\section{Discussion}
\label{S-D}

In this study, we considered a problem of the assessment of 6 criteria used by university students for hotel selection when attending a short-term educational program in a foreign country. The data available for the assessment procedure include the results of a survey of 202 respondents who report on their age, sex and whether they have previously visited the country, and evaluate the criteria by different ways. In the course of evaluation, each respondent directly estimates ratings and indicate ranks of criteria, and then compares the criteria in pairs to provide an input for subsequent derivation of ratings from the pairwise comparisons by an appropriate computational method.

To provide more accuracy of the assessment, we concurrently applied three methods of rating alternatives from pairwise comparisons, where the results of the most common approaches, the method of principal eigenvector and method of geometric mean, were considered in line with results of log-Chebyshev approximation. Moreover, the ratings derived by these methods and corresponding ranks of criteria were compared with the ratings and ranks provided by respondents directly. 

First, we investigated the possible influence of parameters and properties of respondents on the degree of difference and similarity of results provided by different methods of deriving ratings. Specifically, we examined the assumption that the consistency between direct ranks and ratings given by respondents may affect the results of calculating ranks from pairwise comparisons by these methods. As our analysis has shown, for the majority of respondents, the consistency of the direct judgments, as well as parameters of respondents do not have a significant impact on the stability of the results, which made it inappropriate to divide respondents into groups for separate analysis.

Furthermore, we considered ratings provided by respondents directly or indirectly through pairwise comparisons of criteria. Calculation of the mean vectors of ratings has indicated that on average all indirect methods produce very close results. Specifically, the mean vectors of the methods of principal eigenvector and geometric mean appear to be closest to each other, whereas the mean vectors of the best and worst differentiating results of log-Chebyshev approximation form lower and upper bounds for the other two mean vectors. The mean vector of direct ratings by respondents is rather different from the mean vectors of indirect methods, which is due to a rough 5-point scale that narrows the range of assessment values used by respondents.

On average all methods equally determine the highest ratings of criteria in the order $(2,1,5)$. The other criteria are ordered by all indirect methods as $(3,4,6)$, and by the ratings directly provided by respondents as $(6,3,4)$. The standard deviation of ratings of criterion 2 has the smallest or the next smallest value among the criteria, which speaks in favor of a unanimous assessment of this criterion by all respondents. As a result, the analysis of mean vectors of ratings suggests two possible candidates for the optimal rank vector: $(2,1,5,3,4,6)$ and $(2,1,5,6,3,4)$.

To give further insight into the consistency or inconsistency of ratings derived by different methods, we estimated Pearson correlation coefficients between vectors of ratings obtained for all respondents. We have found a high correlation close to one between the results of the methods of principal eigenvector and geometric mean. The correlation between ratings from these two methods and from the log-Chebyshev approximation is somewhat lower, but still close to one. The good agreement of the ratings obtained from pairwise comparisons by all methods can be seen as a result of the natural and meaningful character of the input data that reflect mechanisms of human judgment and hence have a certain level of immanent consistency.

As the next step of the investigation, we examined the vectors of ranks provided by direct and indirect assessments. First, we checked whether the rank vectors $(2,1,5,3,4,6)$ or $(2,1,5,6,3,4)$ derived from the mean vectors of ratings actually occur in the results based on the data provided by individual respondents. We have found that these vectors of ranks are relatively rare as results of all assessment methods, and thus cannot be considered as appropriate solutions. At the same time, calculating the number of the most frequently obtained rank vectors has shown that the most occurred vector is $(2,1,5,3,6,4)$ for all methods except for the ranking from direct ratings by respondents, where this vector is the second most frequent. We note that the last vector of ranks provides the same order for three most important criteria as the two vectors considered above, and differ from them only for the last three least important criteria. 

It follows from the analysis of the results obtained that for the given survey data, one can take the rank vector $(2,1,5,3,6,4)$ as the most accurate assessment of criteria for hotel selection. Assuming the respondents are representative of the total population of university students, we can conclude that the typical student ranks the criteria of hotel selection in the following order (from the most important to the least important): 
\begin{enumerate}
\item
Accommodation cost, $\mathcal{C}_{2}$; 
\item
Hotel location, $\mathcal{C}_{1}$;
\item
Room amenities, $\mathcal{C}_{5}$;
\item
Typical guests, $\mathcal{C}_{3}$;
\item
Courtesy of staff, $\mathcal{C}_{6}$;
\item
Free breakfast, $\mathcal{C}_{4}$.
\end{enumerate}

As another less definite but more realistic conclusion, we can group criteria into two levels. The first level consists of the criterion $\mathcal{C}_{2}$ (accommodation cost), which is ranked first most frequently, and then criteria $\mathcal{C}_{1}$ (hotel location) and $\mathcal{C}_{5}$ (room amenities), which can sometimes change the order of each other. On the second level of lower importance are the criteria $\mathcal{C}_{3}$ (typical guests), $\mathcal{C}_{6}$ (courtesy of staff) and $\mathcal{C}_{4}$ (free breakfast), which can go in almost any order. 

Finally, we estimated Kendall correlation coefficients for rank vectors derived from pairwise comparisons. As in the case of ratings, the correlations between rank vectors have been found close to one, which shows once again that for the pairwise comparison data provided by respondents, the method of principal eigenvector, the method of geometric mean and the method of log-Chebyshev approximation produce quite consistent outcome. The similarity of assessment results can serve as an additional argument in support of accepting the obtained order of criteria as the optimal solution.

\section{Conclusions}
\label{S-C}

The problem of evaluating preferences of criteria for choice is a key component in the multicriteria decision making procedure of assessment alternatives, and it is of independent interest in many applications including marketing research and practice. However, the existing methods to handle the problem, which are based on both direct evaluation of each criterion and indirect evaluation, where the criteria are compared in pairs, are known to give different and even opposite results. In the case when different methods may lead to inconsistent outcomes, a natural approach is to solve the problem by using several most widely used or mathematically justified methods concurrently and then combine the results to provide a more consistent and justified final solution.

In the paper, we have presented an approach that is intended to improve quality of the evaluation of criteria used by students in hotel selection when attending an educational program abroad. The approach is based on implementation of both direct and indirect methods of evaluating criteria, followed by statistical procedures to examine similarity and difference of the results. In addition, we have compared the solutions obtained by three methods of evaluating criteria from pairwise comparisons, including the methods of principal eigenvector, geometric mean and log-Chebyshev approximation. 

By applying this technique, we have derived a solution that ranks the criteria of hotel selection in a way, which combines the results of different methods to provide the basis for more accurate assessment of the criteria. The comparison of the solutions produced by the indirect methods under study has shown a fairly high degree of similarity of results. We consider that the combined assessment technique described in the paper may serve as a template for a useful solution approach to handle real world problems of evaluating alternatives, where the accuracy and reliability of results are of prime importance and thus require additional formal and experimental justification.

As directions of future research, one can consider an extension of the approach to incorporate additional methods of evaluating alternatives from pairwise comparisons, including interval, linguistic and fuzzy pairwise comparisons. Further development of the technique into efficient algorithms and procedures to support decision making processes in practical problems presents another promising line of research.

\bibliographystyle{abbrvurl}

\bibliography{Using_pairwise_comparisons_to_evaluate_consumer_preferences_in_hotel_selection}

\end{document}